\documentclass{article}[12pt]

\usepackage{latexsym}
\usepackage{amsmath}
\usepackage{amssymb}

\begin{document}

\title{A solvable nonlinear autonomous recursion of arbitrary order}

\author{Francesco Calogero$^{a,b}$\thanks{e-mail: francesco.calogero@roma1.infn.it}
\thanks{e-mail: francesco.calogero@uniroma1.it}
 , Farrin Payandeh$^c$\thanks{e-mail: farrinpayandeh@yahoo.com}
 \thanks{e-mail: f$\_$payandeh@pnu.ac.ir}}

\maketitle   \centerline{\it $^{a}$Physics Department, University of
Rome "La Sapienza", Rome, Italy}

\maketitle   \centerline{\it $^{b}$INFN, Sezione di Roma 1}

\maketitle

\maketitle   \centerline{\it $^{c}$Department of Physics, Payame
Noor University, PO BOX 19395-3697 Tehran, Iran}

\maketitle

\begin{abstract}

The \textit{initial-values} problem of the following \textit{nonlinear}
\textit{autonomous} recursion of order $p$,%
\begin{equation*}
z\left( s+p\right) =c\prod\limits_{\ell =0}^{p-1}\left\{ \left[ z\left(
s+\ell \right) \right] ^{a_{\ell }}\right\} ~,~~~
\end{equation*}%
---with $p$ an \textit{arbitrary }positive integer, $z\left( s\right) $ the
dependent variable (possibly a \textit{complex} number), $s$ the independent
variable (a \textit{nonnegative integer}), $c$ an \textit{arbitrarily
assigned}, possibly\textit{\ complex}, number, and the $p$ exponents $%
a_{\ell }$ \textit{arbitrarily assigned integers} (positive, negative or
vanishing, so that the right-hand side of the recursion be \textit{univalent}%
)---is \textit{solvable} by \textit{algebraic} operations, involving the
solution of a system of \textit{linear algebraic} equations (generally
\textit{explicitly solvable}) and of a single \textit{polynomial} equation
of degree $p$ (hence \textit{explicitly solvable} for $p=1,2,3,4$).

\end{abstract}

\section{Introduction}

In this paper we point out that the \textit{initial-values} problem of the
following \textit{nonlinear} \textit{recursion} of order $p$,
\begin{subequations}
\begin{equation}
z\left( s+p\right) =c\prod\limits_{\ell =0}^{p-1}\left\{ \left[ z\left(
s+\ell \right) \right] ^{a_{\ell }}\right\} ~,~~~  \label{1Rec}
\end{equation}%
is \textit{solvable} by \textit{algebraic} operations, involving the
solution of an (generally \textit{explicitly solvable}) system of \textit{%
linear algebraic} equations, and of a single \textit{polynomial} equation of
degree $p$ (hence \textit{explicitly solvable} for $p=1,2,3,4$).
And---trivial as this clarification might be---let us emphasize that
"solving " the initial-values problem of the recursion (\ref{1Rec}) implies
the determination of the function $F\left( z\left( 0\right) ,z\left(
1\right) ,...,z\left( p-1\right) ;c;a_{0},a_{1},...,a_{p-1};s\right) $ such
that, for \textit{all} values of the independent variable $s,$
\begin{equation}
z\left( s\right) =F\left( z\left( 0\right) ,z\left( 1\right) ,...,z\left(
p-1\right) ;c;a_{0},a_{1},...,a_{p-1};s\right) ~,~\ ~s=0,1,2,...,
\end{equation}%
as a consequence of (\ref{1Rec}); \textit{not} solving step-by-step
(obviously each being an easy algebraic operation) the recursion (\ref{1Rec}%
) starting from the assigned $p$ initial values $z\left( 0\right) ,z\left(
1\right) ,...,z\left( p-1\right) $---which is of course \textit{not} the way
to find the solution $z\left( s\right) $ for \textit{all} values of the
independent variable $s,$ this procedure becoming indeed unmanageable for
large $s$.

\textbf{Notation 1-1}. Above and hereafter $z\left( s\right) $ is the
\textit{dependent} variable (generally a \textit{complex} number), $s$ is
the \textit{independent} variable (a \textit{nonnegative integer}), $c$ is
an \textit{arbitrarily assigned}, possibly\textit{\ complex}, number, the $p$
exponents $a_{\ell }$ in the right-hand side of (\ref{1Rec}) are \textit{%
arbitrarily assigned integers} (positive, negative or vanishing, so that the
right-hand side of the recursion be \textit{univalent}; although the
relaxation of this constraint might be worthy of future investigation), and
of course the $p$ (possibly \textit{complex}) numbers $z\left( 0\right)
,z\left( 1\right) ,...,z\left( p-1\right) $ are the \textit{arbitrarily
assigned initial data}.$\ $Below the Kronecker symbol $\delta _{\ell s}$ has
the standard meaning: $\delta _{\ell s}=1$ for $\ell =s,$ $\delta _{\ell
s}=0 $ for $\ell \neq s$.$~\blacksquare $

This paper can be considered a follow-up to the recent papers \cite{CP2021b}-%
\cite{CP2021d}, \cite{PZ2021} identifying solvable Ordinary Differential
Equations (ODEs) and Difference Equations (DEs, or equivalently \textit{%
recursions}); its main novelty is to identify a simple class of \textit{%
nonlinear recursions} of \textit{arbitrary order} which are \textit{solvable}%
---in the sense specified above and below---without requiring any
restriction on their parameters nor on their initial data. The issue of
course is whether---or to what extent---this finding is indeed \textit{new},
due to the difficulty to fully explore the relevant literature (see, for
instance, \cite{B1880}, \cite{C2005}-\cite{M}). We did not find it
explicitly reported in the excellent open-access website \textbf{EqWorld}
managed by A. D. Polyanin, where it certainly belongs. But see below \textbf{%
Remark 3-1}.

In the following \textbf{Section 2} our main results are reported; the
corresponding proofs---to the extent they are needed---are provided in
\textbf{Appendix A}. And a quite terse \textbf{Section 3} mentions possible
future developments.

\section{Results}

In this Section we report the main results of this paper; in the interest of
pedagogical clarity, we firstly consider the simplest case with $p=1,$ then
with $p=2,$ and finally we exhibit the result for \textit{arbitrary positive
integer} $p,$ where $p$ is of course the order of the recursion.

\subsection{p=1}

\textbf{Proposition 2-1}. The solution of the \textit{initial-value} problem
for the recursion
\end{subequations}
\begin{equation}
z\left( s+1\right) =c\left[ z\left( s\right) \right] ^{a}~~~~
\label{Prop21Rec}
\end{equation}%
reads as follows:
\begin{subequations}
\label{Prop21Sol}
\begin{equation}
z\left( s\right) =c^{\left( 1-a^{s}\right) /\left( 1-a\right) }\left[
z\left( 0\right) \right] ^{a^{s}}~;  \label{Prop21Sola}
\end{equation}%
unless
\begin{equation}
a=1,
\end{equation}%
in which case the solution reads as follows:%
\begin{equation}
z\left( s\right) =c^{s}\left[ z\left( 0\right) \right] ^{s+1}~.~\blacksquare
\end{equation}

This result is simple enough to make the verification that indeed the
formulas (\ref{Prop21Sol}) provides the solution of the \textit{initial-value%
} problem of the recursion (\ref{Prop21Rec}) a quite simple exercise; we
nevertheless provide an explanation of how it has been obtained below, see
\textbf{Appendix A}.

\textbf{Remark 2.1-1}. It is quite obvious that the result described by
\textbf{Proposition 2-1} holds just as well if (only) the dependent variable
$z\left( s\right) $, instead of being a \textit{scalar} entity, were, for
instance, a \textit{matrix} of \textit{arbitrary} rank; thereby involving in
the recursion (\ref{Prop21Rec}) a much larger number of (scalar) dependent
variables, i. e. \textit{all} the entries of that matrix. $\blacksquare $

\subsection{p=2}

\textbf{Proposition 2-2}. The solution of the \textit{initial-values}
problem for the recursion
\end{subequations}
\begin{equation}
z\left( s+2\right) =c~\left[ z\left( s+1\right) \right] ^{a_{1}}~\left[
z\left( s\right) \right] ^{a_{0}}  \label{p21Rec}
\end{equation}%
is given---for \textit{generic} values of the $2$ (integer) parameters $%
a_{0} $ and $a_{1}$; for some \textit{special} values see below---by the
following formulas:
\begin{subequations}
\label{Prop22}
\begin{equation}
z\left( s\right) =c^{\gamma \left( s\right) }~\left[ z\left( 1\right) \right]
^{\alpha _{1}\left( s\right) }~\left[ z\left( 0\right) \right] ^{\alpha
_{0}\left( s\right) }~,  \label{Prop1zs}
\end{equation}%
with%
\begin{equation}
\alpha _{1}\left( s\right) =\left[ \left( y_{+}\right) ^{s}-\left(
y_{-}\right) ^{s}\right] /d~,  \label{Prop1alf1s}
\end{equation}%
\begin{equation}
\alpha _{0}\left( s\right) =-y_{+}y_{-}~\left[ \left( y_{+}\right)
^{s-1}-\left( y_{-}\right) ^{s-1}\right] /d,  \label{Prop1alf0s}
\end{equation}%
\begin{eqnarray}
\gamma \left( s\right) &=&\left[ \left( 1-a_{1}-a_{0}\right) d\right]
^{-1}\left\{ d+\left( y_{-}-1\right) \left( y_{+}\right) ^{s}-\left(
y_{+}-1\right) \left( y_{-}\right) ^{s}\right\}  \nonumber \\
&=&\left( 1-a_{1}-a_{0}\right) ^{-1}\left[ 1-\alpha _{1}\left( s\right)
-\alpha _{0}\left( s\right) \right] ~,  \label{Prop1gamma}
\end{eqnarray}%
where%
\begin{equation}
y_{\pm }=\left( a_{1}\pm d\right) /2~,~~~d=\sqrt{\left( a_{1}\right)
^{2}+4a_{0}}~.  \label{21y+-d}
\end{equation}

This finding remains true also in the special case with
\end{subequations}
\begin{subequations}
\label{a0Plusa1^2Eq0}
\begin{equation}
a_{0}=-\left( a_{1}/2\right) ^{2}~,
\end{equation}%
when $d$ vanishes, $d=0,$ by taking the appropriate limit as $d\rightarrow 0$%
. Then the solution of the recursion (\ref{p21Rec}) is still given by the
expression (\ref{Prop1zs}), but now with
\begin{equation}
\alpha _{1}\left( s\right) =sy^{s-1}~,~~~\alpha _{0}\left( s\right) =\left(
1-s\right) y^{s}~,~~~\gamma \left( s\right) =\frac{1-\alpha _{1}\left(
s\right) -\alpha _{0}\left( s\right) }{1-a_{1}-a_{0}}~,  \label{al1&al2&gam}
\end{equation}%
where%
\begin{equation}
y=a_{1}/2~.  \label{gam}
\end{equation}

While, in the special case
\end{subequations}
\begin{subequations}
\begin{equation}
a_{1}\neq 2~,~~~a_{1}+a_{0}=1~,  \label{Prop22a0=1-a1}
\end{equation}%
the solution of the initial-values problem of the recursion (\ref{p21Rec})
reads as follows:%
\begin{equation}
z\left( s\right) =z\left( 0\right) ~c^{\gamma \left( s\right) }~\left[
z\left( 1\right) /z\left( 0\right) \right] ^{\left[ 1-\left( a_{1}-1\right)
^{s}\right] /\left( 2-a_{1}\right) }~,  \label{zs}
\end{equation}%
with%
\begin{equation}
\gamma \left( 0\right) =\gamma \left( 1\right) =0~;~~~\gamma \left( s\right)
=\sum_{\ell =0}^{s-2}\left[ \left( s-1-\ell \right) \left( a_{1}-1\right)
^{\ell }\right] ~,~s>2~.  \label{Prop1gam}
\end{equation}

Finally, in the even more special case with
\end{subequations}
\begin{subequations}
\label{a1Eq2&a0Eq-1}
\begin{equation}
a_{1}=2~,~~~a_{0}=-1~,
\end{equation}%
the solution of the initial-values problem of the recursion (\ref{p21Rec})
reads as follows:
\begin{equation}
z\left( s\right) =c^{s\left( s-1\right) /2}~\left[ z\left( 1\right) \right]
^{s}~\left[ z\left( 0\right) \right] ^{1-s}~.~\blacksquare
\end{equation}

For an indication of how these findings have been arrived at see \textbf{%
Appendix A}; but to check their validity it is sufficient to verify that in
each case the indicated \textit{solution} does indeed satisfy the relevant
\textit{recursion} with the \textit{arbitrarily assigned initial data }$%
z_{1}\left( 0\right) ,$ $z_{2}\left( 0\right) $.

\subsection{p is an \textit{arbitrary positive integer}}

The results reported below are for \textit{generic} values of all involved
parameters. And their validity shall be evident---after the explanations
provided below---to all readers who have previously digested the results
reported in the $2$ preceding Subsections.

\textbf{Proposition 2-3}.\textbf{\ }The solution of the \textit{%
initial-values} problem of the recursion (\ref{1Rec}) is provided by the
following formula
\end{subequations}
\begin{equation}
z\left( s\right) =c^{\gamma \left( s\right) }~\prod\limits_{\ell
=0}^{p-1}\left\{ \left[ z\left( \ell \right) \right] ^{\alpha _{\ell }\left(
s\right) }\right\} ~,  \label{23Sol}
\end{equation}%
where of course the $p$ values $z\left( \ell \right) $ with $\ell
=0,1,...,p-1$ are the\textit{\ arbitrarily assigned} \textit{initial data},
while the function $\gamma \left( s\right) ,$ as well as the $p$ functions $%
\alpha _{\ell }\left( s\right) $ (with $\ell =0,1,...,p-1$), can be
calculated as shown below.

But firstly we note that this formula (\ref{23Sol}) clearly implies that the
exponents $\gamma \left( s\right) $ and $\alpha _{\ell }\left( s\right) $
satisfy the following \textit{initial conditions}:
\begin{subequations}
\label{23IniAlfGam}
\begin{equation}
\gamma \left( s\right) =0~,~~~s=0,1,...,p-1~,
\end{equation}%
\begin{equation}
\alpha _{n}\left( s\right) =\delta _{ns}~,~~~n,s=0,1,...,p-1~.
\end{equation}

Next, we define the $p$ quantities $y_{\ell }$, with $\ell =0,1,...,p-1,$ as
the $p$ roots of the algebraic equation
\end{subequations}
\begin{equation}
y^{p}-\sum_{\ell =0}^{p-1}\left( a_{\ell }~y^{\ell }\right) =0~;
\end{equation}%
and we hereafter assume that they are \textit{all different}---consistently
with our \textit{genericity} premise, introduced to avoid the consideration
of a plethora of special subcases. Of course for $p=1,2,3,4$ these roots $%
y_{\ell }$ can be computed \textit{explicitly} in terms of the $p$
parameters $a_{\ell }$ (and, in some special cases, also for $p=6$ and $p=8$
or $p=9$; see \cite{CPSextic} and \cite{CPEqDeg8and9}).

The exponents $\gamma \left( s\right) ,$ as well as the $p$ exponents $%
\alpha _{\ell }\left( s\right) ,$ featured by the right-hand side of the
solution formula (\ref{23Sol}), are then clearly the solutions of the
following \textit{linear recursions}:
\begin{subequations}
\begin{equation}
\alpha _{n}\left( s+p\right) =\sum_{\ell =0}^{p-1}\left[ a_{\ell}~\alpha
_{n}\left( s+\ell \right) \right] ~,~~~n=0,1,2,...,p-1~,
\end{equation}%
\begin{equation}
\gamma \left( s+p\right) =1+\sum_{\ell =0}^{p-1}\left[ a_{\ell}~\gamma
\left( s+\ell \right) \right] ~~;
\end{equation}%
hence they are given by the following \textit{explicit} formulas:
\end{subequations}
\begin{equation}
\alpha _{n}\left( s\right) =\sum_{\ell =0}^{p-1}\left[ A_{n\ell }~\left(
y_{\ell }\right) ^{s}\right] ~,~~~n=0,1,...,p-1~,
\end{equation}%
\begin{equation}
\gamma \left( s\right) =\left[ 1-\sum_{\ell =0}^{p-1}\left( a_{\ell }\right) %
\right] ^{-1}+\sum_{\ell =0}^{p-1}\left[ C_{\ell }~\left( y_{\ell }\right)
^{s}\right] ~,
\end{equation}%
where the $p^{2}$ ($s$-independent) parameters $A_{nl}$ and the $p$ ($s$%
-independent) parameters $C_{\ell }$ are defined by the following system of
\textit{linear} equations implied by the \textit{initial values} (\ref%
{23IniAlfGam}):
\begin{subequations}
\begin{equation}
\sum_{\ell =0}^{p-1}\left[ A_{n\ell }~\left( y_{\ell }\right) ^{s}\right]
~=\delta _{ns}~,~~n,s=0,1,2,...,p-1~,~
\end{equation}%
\begin{equation}
\sum_{\ell =0}^{p-1}\left[ C_{\ell }~\left( y_{\ell }\right) ^{s}\right] =-%
\left[ 1-\sum_{\ell l=0}^{p-1}\left( a_{\ell }\right) \right] ^{-1}~.
\end{equation}%
These \textit{linear} algebraic equations are of course \textit{solvable},
providing \textit{explicit}---if, for large $p$, quite
complicated---expressions of the $p^{2}+p$ parameters $A_{n\ell }$ and $%
C_{\ell }$. $\blacksquare $

\textbf{Remark 2-1}. A rather trivial extension of the findings reported
above is via the change of dependent variables
\end{subequations}
\begin{equation}
z\left( s\right) =\tilde{z}\left( s\right) +f\left( s\right) ~,
\end{equation}%
where $\tilde{z}\left( s\right) $ is the new dependent variable and $f\left(
s\right) $ an \textit{arbitrarily} chosen function; or just an $s$%
-independent parameter, in order to maintain the \textit{autonomous}
character of the resulting new recursion satisfied by $\tilde{z}\left(
s\right) $. $\blacksquare $

\section{Outlook}

The generalization of the findings reported in this paper to \textit{%
nonautonomous} recursions (by other means than that mentioned in \textbf{%
Remark 2-1}), as well as to \textit{systems} of nonlinear recursions
involving \textit{several} dependent variables, are interesting tasks; and
in that context the possibility to also consider the case of \textit{%
noncommuting} dependent variables (such as matrices: see, for instance, \cite%
{MikSok} \cite{VS0}) is another interesting possible development (it is
already available in the $p=1$ case: see \textbf{Remark 2.1-1}).

\textbf{Remark 3-1}. This important remark has been added after the present
paper was completed. The recursion (\ref{1Rec}) which we investigate in this
paper is evidently \textit{nonlinear}, but---as pointed out to us by Paolo
Santini---it is easily transformed into a \textit{linear} recursion via the
following simple change of dependent variables:
\begin{subequations}
\begin{equation}
\ln \left[ z\left( s\right) \right] =x\left( s\right) ~,  \label{lnzsToalf}
\end{equation}%
which indeed transforms the \textit{nonlinear} recursion (\ref{1Rec})
satisfied by the dependent variable $z\left( s\right) $ into the following
\textit{linear }recursion satisfied by the new dependent variable $x\left(
s\right) $:%
\begin{equation}
x\left( s+p\right) =\ln \left( c\right) +\sum_{\ell =0}^{p-1}\left[
a_{\ell}~x\left( s+\ell \right) \right] ~\ .  \label{alfaRec}
\end{equation}%
It is indeed this kind of change of variables that subtends the findings
reported in the present paper, see in particular the \textit{ansatz} (\ref%
{23Sol}). Note however that the solutions of these \textit{linear}
recursions (\ref{alfaRec}) entail themselves the solution of \textit{%
nonlinear algebraic} operations, as indeed displayed by the solutions
reported above.

Finally let us emphasize---to complement what was stated above (in the
next-to-last paragraph of \textbf{Section 1})---that the transition from (%
\ref{1Rec}) to (\ref{alfaRec}) via (\ref{lnzsToalf}) implies that the
solvability of (\ref{1Rec}) can be reduced to the solvability of (\ref%
{alfaRec}); which is indeed mentioned in \textbf{EqWorld}, see there \textit{%
Exact Solutions $>$ Functional Equations $>$ 1. Linear
Difference and Functional Equations with One Independent Variable
$>$ 1.2. Other Linear Difference and Functional Equations
$>$ see Items 16 and 15}; although a few more developments are
still needed to go from the general solution described in \textbf{EqWorld}
to the solution of the \textit{initial-values} problem discussed in our paper%
\textit{. }$\blacksquare $

\section{Appendix A}

In this \textbf{Appendix} we explain how the results---\textbf{Propositions
2-1} and \textbf{2-2}---reported in the $2$ \textbf{Subsections 2.1} and
\textbf{2.2} have been arrived at. The following subdivision in Subsections
indeed corresponds to the analogous subdivision in \textbf{Section 2}. As
for the main result, as reported in \textbf{Subsection 2.3, }we believe that
a more detailed explanation than that provided in the formulation of \textbf{%
Proposition 2-3} is unnecessary, given the analogy of that treatment to that
provided below in the case with $p=2$.

\subsection{p=1}

In this Subsection we consider the simple recursion
\end{subequations}
\begin{equation}
z\left( s+1\right) =c\left[ z\left( s\right) \right] ^{a}  \label{pEq1Rec}
\end{equation}%
(in the case with $a\neq 1$; the case with $a=1$ being too trivial to
require an explicit treatment).

We then introduce the following \textit{ansatz}
\begin{subequations}
\begin{equation}
z\left( s\right) =c^{\gamma \left( s\right) }\left[ z\left( 0\right) \right]
^{\alpha \left( s\right) }  \label{21Ans}
\end{equation}%
clearly implying%
\begin{equation}
\gamma \left( 0\right) =0~,~~~\alpha \left( 0\right) =1~.  \label{gamalf0}
\end{equation}

We then insert this \textit{ansatz} (\ref{21Ans}) in the recursion (\ref%
{pEq1Rec}) and thereby easily get the following \textit{linear} relations:
\begin{equation}
\gamma \left( s+1\right) =1+a\gamma \left( s\right)  \label{21gamasa}
\end{equation}%
implying, via the position%
\begin{equation}
\gamma \left( s\right) =\tilde{\gamma}\left( s\right) +1/\left( 1-a\right) ~,
\label{gamtilde}
\end{equation}%
with (via (\ref{gamalf0}))%
\begin{equation}
\tilde{\gamma}\left( 0\right) =-1/\left( 1-a\right) ~,  \label{gamtilde0}
\end{equation}%
the recursion%
\begin{equation}
\tilde{\gamma}\left( s+1\right) =a~\tilde{\gamma}\left( s\right) ~;
\end{equation}%
as well as
\begin{equation}
\alpha \left( s+1\right) =a~\alpha \left( s\right) ~.
\end{equation}%
The last $2$ recursions clearly imply%
\begin{equation}
\alpha \left( s\right) =\alpha \left( 0\right) ~a^{s}~,~~~\tilde{\gamma}%
\left( s\right) =\tilde{\gamma}\left( 0\right) ~a^{s}\
\end{equation}%
hence, via (\ref{gamtilde}) and (\ref{gamtilde0}),%
\begin{equation}
\gamma \left( s\right) =\left( 1-a^{s}\right) /\left( 1-a\right) ~.
\end{equation}%
The results reported in \textbf{Proposition 2-1} are thus proven.

\subsection{p=2}

In this Subsection we show how the findings reported in \textbf{Proposition
2-2} were obtained. Our focus is of course on the recursion (\ref{p21Rec}).

We then introduce the following \textit{ansatz}:
\end{subequations}
\begin{subequations}
\label{p2Ans}
\begin{equation}
z\left( s\right) =c^{\gamma \left( s\right) }\left[ z\left( 1\right) \right]
^{\alpha _{1}\left( s\right) }\left[ z\left( 0\right) \right] ^{\alpha
_{0}\left( s\right) }  \label{p2Ansz}
\end{equation}%
clearly implying%
\begin{eqnarray}
\gamma \left( 0\right) &=&0~,~~~\alpha _{1}\left( 0\right) =0~,~~~\alpha
_{0}\left( 0\right) =1~,  \nonumber \\
\gamma \left( 1\right) &=&0~,~~~\alpha _{1}\left( 1\right) =1~,~~~\alpha
_{0}\left( 1\right) =0~.  \label{p2Ansb}
\end{eqnarray}

The introduction of this \textit{ansatz} in the recursion (\ref{p21Rec})
yields the following simple \textit{linear} \textit{decoupled} recursions:
\end{subequations}
\begin{subequations}
\label{21RecGamAlf}
\begin{equation}
\gamma \left( s+2\right) =1+a_{1}\gamma \left( s+1\right) +a_{0}\gamma
\left( s\right) ~,  \label{21RecGam}
\end{equation}%
\begin{equation}
\alpha _{n}\left( s+2\right) =a_{1}\alpha _{n}\left( s+1\right) +a_{0}\alpha
_{n}\left( s\right) ~,~~~n=1,0~.  \label{21RecAlf}
\end{equation}

The solutions of the $2$ recursions (\ref{21RecAlf}) clearly read as
follows:
\end{subequations}
\begin{subequations}
\begin{equation}
\alpha _{n}\left( s\right) =A_{n+}\left( y_{+}\right) ^{s}+A_{n-}\left(
y_{-}\right) ^{s}~~,~~~n=1,0~,
\end{equation}%
with $y_{\pm }$ the $2$ solutions of the quadratic equation~%
\begin{equation}
y^{2}-a_{1}y-a_{0}=0
\end{equation}%
implying the expression (\ref{21y+-d}) of $y_{\pm }$ and $d$. In these
formulas the $4$ parameters $A_{n\pm }$ (with $n=1,0$) are \textit{a priori
arbitrary}, but they are clearly related to the $4$ initial values $\alpha
_{n}\left( 0\right) ,$ $\alpha _{n}\left( 1\right) $ (again, with $n=1,0$)
by the following $4$ formulas:
\end{subequations}
\begin{subequations}
\begin{equation}
A_{n+}+A_{n-}=\alpha _{n}\left( 0\right) ~,~~~A_{n+}y_{+}+A_{n-}y_{-}=\alpha
_{n}\left( 1\right) ~,~~~n=0,1~,
\end{equation}%
implying, via (\ref{p2Ansb}),%
\begin{eqnarray}
A_{1+}+A_{1-} &=&0~,~~~A_{1+}y_{+}+A_{1-}y_{-}=1~,  \nonumber \\
A_{0+}+A_{0-} &=&1~,~~~A_{0+}y_{+}+A_{0-}y_{-}=0~,
\end{eqnarray}%
hence%
\begin{equation}
A_{1\pm }=\pm 1/d~,~~~A_{0\pm }=\mp y_{\mp }/d~,
\end{equation}%
of course with $d=\sqrt{\left( a_{1}\right) ^{2}+4a_{0}}$ (see (\ref{21y+-d}%
)). Which clearly imply the $2$ relations (\ref{Prop1alf1s}) and (\ref%
{Prop1alf0s}).

The next step is to determine the function $\gamma \left( s\right) .$ The
position
\end{subequations}
\begin{equation}
\gamma \left( s\right) =\tilde{\gamma}\left( s\right) +1/\left(
1-a_{0}-a_{1}\right) ~,
\end{equation}%
when inserted in the recursion (\ref{21RecGam}) satisfied by $\gamma \left(
s\right) ,$ implies that $\tilde{\gamma}\left( s\right) $ satisfies the
\textit{same} recursion (\ref{21RecAlf}) satisfied by $\alpha _{n}\left(
s\right) $; hence a repetition of the development detailed just above leads
to the solution (\ref{Prop1gam}) for $\gamma \left( s\right) ,$ thereby
completing the derivation of the first part of \textbf{Proposition 2-2}.

As for the remaining part of \textbf{Proposition 2-2}, it is easily seen
that the formulas (\ref{al1&al2&gam}) follow from the formulas (\ref{Prop22}%
) by taking the limit $a_{0}\rightarrow -\left( a_{1}/2\right) ^{2}$
implying $d\rightarrow 0,$ and likewise the formulas (\ref{zs}) with (\ref%
{Prop1gam}) follow from the formulas (\ref{Prop22}) by taking the limit $%
a_{0}\rightarrow 1-a_{1}$ implying $d\rightarrow a_{1}-2,$ $y_{+}\rightarrow
a_{1}-1,$ $y_{-}\rightarrow 1;$ while a check of the very special case (\ref%
{a1Eq2&a0Eq-1}) is quite trivial.

\textbf{Acknowledgements}. We like to acknowledge with many thanks the
insight by prof. Paolo Santini mentioned in \textbf{Remark 3-1}, as well as
useful discussions with prof. Andrea Giansanti on the application of
recursions. And we also like to acknowledge with thanks $3$ grants,
facilitating our collaboration---mainly developed via e-mail exchanges---by
making it possible for FP to visit $3$ times (one of which in the future)
the Department of Physics of the University of Rome "La Sapienza": $2$
granted by that University, and one granted jointly by the Istituto
Nazionale di Alta Matematica (INdAM) of that University and by the
International Institute of Theoretical Physics (ICTP) in Trieste in the
framework of the ICTP-INdAM "Research in Pairs" Programme.\textbf{\ }%
Finally, we also\ like to thank Fernanda Lupinacci who, in these difficult
times---with extreme efficiency and kindness---facilitated all the
arrangements necessary for the presence of FP with her family in Rome.

\end{document}